
\baselineskip=14pt
\parskip=10pt

\magnification=\magstephalf
\def\T{{\cal T}}

\def\1{{\overline{1}}}
\def\2{{\overline{2}}}
\parindent=0pt
\overfullrule=0in

\def\frac#1#2{{#1 \over #2}}

\centerline
{
 \bf 
On The Limiting Distributions of the Total Height On Families of Trees
}

\bigskip
\centerline
{\it By Andrew LOHR and Doron ZEILBERGER}

{\bf Abstract}: A symbolic-computational algorithm, fully implemented in Maple, is described, that computes
explicit expressions for generating functions that enable the efficient computations of the expectation, variance, and higher moments, of
the random variable `sum of distances to the root', defined on any given family of rooted ordered trees (defined
by degree restrictions). 
Taking limits, we confirm, via elementary methods, the fact, due to David Aldous, and expanded by Svante Janson and
others, that the limiting (scaled) distributions are all the
same, and coincide with the limiting distribution of the same random variable, when it is defined on labeled rooted trees.

{\bf Maple packages and Sample Input and Output Files}

This article is accompanied by  Maple packages, {\tt TREES.txt}, 
and {\tt THS.txt},
and several input and output files available from
the front of this article

{\tt http://www.math.rutgers.edu/\~{}zeilberg/mamarim/mamarimhtml/otrees.html} \quad .

{\bf Background} 

While many natural families of combinatorial random variables, $X_n$, indexed by a positive integer $n$,
(for example, tossing a coin $n$ times and noting the number of Heads, or counting the number of occurrences of a specific pattern
in an $n$-permutation) have  different expectations, $\mu_n$, and different standard deviations, $\sigma_n$, and (usually) largely different
asymptotic expressions for these, yet the centralized and scaled versions, $Z_n:=\frac{X_n -\mu_n}{\sigma_n}$, very often,
converge (in distribution) to the standard normal distribution whose probability density function is famously
$\frac{1}{\sqrt{2 \pi}} exp(-\frac{x^2}{2})$, and whose moments are $0,1,0,3,0,5,0,15,0,105, \dots$.
Such sequences of random variables are called {\it asymptotically normal}.
Whenever this is {\bf not} the case, it is a cause for excitement. One celebrated case
(see Dan Romik's [Ro] masterpiece for an engaging and detailed description) is the random variable
`largest increasing subsequence', defined on the set of permutations, where the intriguing {\it Tracy-Widom distribution} shows up.

Other, more recent, examples of {\it abnormal} limiting distributions are described in [Z1], [EZ1],[EZ2], and [EZ3].

In this article we consider, from an elementary, explicit, {\it symbolic-computational}, viewpoint, the
random variable `sum of distances to the root', defined over an {\it arbitrary} family of ordered rooted trees
defined by degree restricions.

It turns out that the families of trees considered here are
special cases of Galton-Watson trees. These have been studied extensively by continuous probability theorists for many years,
 with a nice, comprehensive introduction given by Janson in [J3].
In particular, they are trees that are determined by determining the number of childeren that every node has by 
independently sampling some fixed distribution with expected value at most 1. Like the trees
considered here (described below), they are also types of Galton-Watson trees. It was shown in 
[A1], [A2], and [MM] that all Galton-Watson generated from a 
finite variance distribution of vertex degrees followed the same distribution as 
the area under a Brownian excursion, also a topic well studied in advanced probability theory.
In particular, Janson, in section 14 of [J1], presents a complicated infinite sum which converges to this distribution originally discovered by Darling (1983).
Asymptotic analysis of mean, variance, and higher moments for Galton-Watson trees can be found in [J4].

All these authors used {\it continuous}, advanced, probability theory, that while very powerful, only gives you the
limit. We are interested in {\it explicit} expressions for the first few moments themselves, or failing this
for explicit expressions for the generating functions, for {\it any} family of rooted ordered trees given by
degree restrictions. In particular, we study in detail the case of {\it complete binary trees}, famously counted
by the Catalan numbers.

We proceed in the same vein as in [EZ2].
In that article, the random variable `sum of the distances from the root', 
defined on the set of {\it labelled rooted trees} on $n$ vertices,
was considered, and it was shown how to find explicit expressions for any given moment, and the
first $12$ moments were derived, extending the pioneering work of John Riordan and Neil Sloane ([RiS]),
who derived an explicit formula for the expectation. 
The exact and approximate values for the {\it limits}, as $n \rightarrow \infty$, of
$\alpha_3$ (the {\it skewness}),  $\alpha_4$ (the {\it kurtosis}),  and the higher moments
through the ninth turn out to be as follows.

$$
\alpha_3 \, = \,
{\frac{ \left( 6\,\pi -{\frac {75}{4}} \right) \sqrt {3}\sqrt {{\frac {\pi }{10-3\,\pi }}}}{10-3\,\pi }}
\, =\,
 0.7005665293596503\dots \quad , 
$$
$$
\alpha_4 \, = \, 
{\frac{-189\,{\pi }^{2}+315\,\pi +884}{ 7 \, \left( 10-3\,\pi  \right) ^{2}}}
\, = \,
3.560394897132889\dots 
\quad , 
$$
$$
\alpha_5 \, \, = \, 
{\frac { \left( 36\,{\pi }^{2}+{\frac {75}{2}}\,\pi -{\frac {105845}{224}} \right) \sqrt {3}\sqrt {{\frac {\pi }{10-3\,\pi }}}}{ \left( 10-3\,\pi 
 \right) ^{2}}}
\, =\, 
 7.2563753582799571\dots \quad , 
$$
$$
\alpha_6 \, \, = \, 
{\frac{15}{16016}}\,{\frac {-144144\,{\pi }^{3}-720720\,{\pi }^{2}+3013725\,\pi +2120320}{ \left( 10-3\,\pi  \right) ^{3}}}
\, = \, 
 27.685525695770609\dots \quad , 
$$
$$
\alpha_7 \, \, = \, 
{\frac{ \left( 162\,{\pi }^{3}+{\frac {6615}{4}}\,{\pi }^{2}-{\frac {103965}{32}}\,\pi -{\frac {101897475}{9152}} \right) \sqrt {3}\sqrt {{\frac {\pi }{
10-3\,\pi }}}}{ \left( 10-3\,\pi  \right) ^{3}}}
\, = \, 
90.0171829093603301\dots \quad, 
$$
$$
\alpha_8 \, \, = \, 
{\frac{3}{2586584}}\,{\frac {-488864376\,{\pi }^{4}-8147739600\,{\pi }^{3}-455885430\,{\pi }^{2}+86568885375\,\pi +32820007040}{ \left( 10-3\,\pi 
 \right) ^{4}}} 
$$
$$
\, = \, 358.80904151261251\dots \quad ,
$$
$$
\alpha_9 \, \, = \, 
{\frac{ \left( 648\,{\pi }^{4}+15795\,{\pi }^{3}+{\frac {591867}{16}}\,{\pi }^{2}-{\frac {461286225}{2288}}\,\pi -{\frac {188411947088175}{662165504}}
 \right) \sqrt {3}\sqrt {{\frac {\pi }{10-3\,\pi }}}}{ \left( 10-3\,\pi  \right) ^{4}}}
\, = \, 1460.7011342971821\dots \quad .
$$

{\bf Acknowledgement}: Many thanks are due to  Valentin F\'eray and Svante Janson for telling us about the
work of Aldous, Marckert and Mokkadem, and Janson.

{\bf This Article}

In this article we extend the work of [EZ2] and treat infinitely many other families of trees.
For any given set of positive integers, $S$, we will have a `sample space' of all ordered rooted trees where
a vertex may have no children (i.e. be a {\it leaf}) or it {\bf must} have a number of children that belongs to $S$.
If $S=\{2\}$ we have the case of {\it complete binary trees}.

For each such family, defined by $S$, we will show how to derive explicit expressions for the generating functions
of the numerators of the straight moments, from which one can easily get many values, and very efficiently find the numerical
values for the moments-about-the-mean and hence the scaled moments. For the special case of complete binary
trees, we will derive explicit expressions for the first nine moments (that may be extended indefinitely),
as well as explicit expressions for the asymptotics of the scaled moments, 
and surprise! they coincide {\it exactly} with those found in [EZ2]
for the case of labelled rooted trees. This leads us to conjecture that the limiting distribution is the {\it same}
for each such family. 

{\bf Rooted Ordered Trees}

Recall that an {\it ordered rooted tree} is an unlabeled graph with the root drawn at the top, and each
vertex has a certain number (possibly zero) of children, drawn from left to right. For any finite set
of positive integers, $S$, let  $\T(S)$ be the set of all rooted labelled trees where each vertex either has
no children, or else has a number of children that belongs to $S$. The set $\T(S)$ has the
following structure (``grammar'')
$$
\T(S) = \{\cdot\} \bigcup_{i \in S}  \, \{\cdot\} \times \T(S)^i \quad.
$$

Fix $S$, Let $f_n$ be number of rooted ordered trees in $\T(S)$ with exactly $n$ vertices.
It follows immediately, by elementary generatingfunctionology, that
the ordinary generating function
$$
f(x) :=\sum_{n=0}^{\infty} f_n  \, x^n \quad ,
$$
(that is the sum of the weights of {\it all} members of $\T(S)$ with the weight $x^{NumberOfVertices}$ assigned to each tree)
satisfies the {\bf algebraic} equation
$$
f(x) = x \left ( 1+ \sum_{i \in S} f(x)^i \right ) \quad .
$$

Given an ordered tree, $t$, define the random variable $H(t)$ to be the sum of the distances to the root of all vertices.
Let $H_n$ be its restriction to the subset of $\T(S)$, let's call it $\T_n(S)$, of members of $\T(S)$ with exactly
$n$ vertices. Our goal in this article is to describe a symbolic-computational algorithm that, for {\it any}
finite set $S$ of positive integers, {\it automatically} finds  generating functions that enable the fast
computation of the average, variance, and as many higher moments as desired. We will be particularly interested
in the limit, as $n \rightarrow \infty$, of the centralized-scaled distribution, and we have strong evidence
to conjecture that it is always the same as the one for rooted labelled trees found in [EZ2].

Let $P_n(y)$ be the generating polynomial defined over $\T_n(S)$, of the random variable, `sum of distances from the root'.
Define the {\it grand generating function}
$$
F(x,y)=\sum_{n=0}^{\infty} P_n(y) x^n \quad .
$$

Consider a typical tree, $t$, in $\T_n(S)$, and now define the more general {\it weight} by $x^{NumberOfVertices}\, y^{H(t)}=x^n \, y^{H(t)}$.
If $t$ is a singleton, then its weight is simply $x^1 y^0=x$, but if its sub-trees (the trees whose roots are the children
of the original root) are $t_1, t_2, \dots t_i$ (where $i \in S$), then 
$$
H(t)=H(t_1)+ \dots + H(t_i) + n-1 \quad,
$$
since when you make the tree $t$, out of  subtrees $t_1, \dots, t_i$ by placing them from left to right and then
attaching them to the root, each vertex gets its `distance to the root' increased by $1$, so altogether the 
sum of the vertices' heights gets increased
by the total number of vertices in $t_1, \dots, t_i$ (i.e. $n-1$). Hence $F(x,y)$ satisfies the
{\bf functional equation}
$$
F(x,y)=x \cdot  \left ( 1+ \sum_{i \in S} F(xy, y)^i \right ) \quad,
$$
that can be used to generate many terms of the sequence of generating polynomials $\{ P_n(y) \}$.

Note that when $y=1$, $F(x,1)=f(x)$, and we get back the algebraic equation satisfied by $f(x)$.

{\bf From Enumeration to Statistics in General}

Suppose that we have a finite set, $A$, on which a certain numerical attribute, called {\it random variable}, $X$,
(using the probability/statistics lingo), is defined.

For any non-negative integer $i$, let's define
$$
N_i:=\sum_{a \in A} X(a)^i \quad .
$$
In particular, $N_0(X)$ is the number of elements of $A$.

The expectation of $X$, $E[X]$, denoted by $\mu$,  is, of course,
$$
 \mu  \, = \, \frac{N_1}{N_0} \quad .
$$

For $i>1$, the $i$-th straight moment is
$$
E[X^i] \, = \, \frac{N_i}{N_0} \quad .
$$

The $i$-th {\it moment about the mean} is
$$
m_i:=E[(X-\mu)^i]= E[\sum_{r=0}^{i} {{i} \choose {r}} (-1)^r \mu^r X^{i-r}]=
\sum_{r=0}^{i}  (-1)^r {{i} \choose {r}} \mu^r E[X^{i-r}]
$$
$$
=\, \sum_{r=0}^{i}  (-1)^r {{i} \choose {r}} \left ( \frac{N_1}{N_0} \right )^r   \frac{N_{i-r}}{N_0}
$$
$$
= \, \frac{1}{N_0^i} \sum_{r=0}^{i}  (-1)^r {{i} \choose {r}} N_1^r N_0^{i-r-1} N_{i-r}  \quad .
$$

{\bf Finally}, the most interesting quantities, statistically speaking, apart from the mean $\mu$ and variance $m_2$ are
the {\bf scaled-moments}, also known as, {\it alpha coefficients}, defined by
$$
\alpha_i :=\frac{m_i}{m_2^{i/2}} \quad .
$$

{\bf Using Generating functions}

In our case $X$ is $H_n$ (the sum of the vertices' distances to the root, defined over rooted ordered trees in our family,
with $n$ vertices), and we have
$$
N_1(n) = P_n'(1)
$$
$$
N_i(n) = (y\frac{d}{dy})^i P_n(y) \bigl \vert_{y=1} .
$$
It is more convenient to first find the numerators of the factorial moments
$$
F_i(n)=(\frac{d}{dy})^i P_n(y) \vert_{y=1} \quad,
$$
from which $N_i(n)$ can be easily found, using the Stirling numbers of the second kind.

{\bf Automatic Generation of Generating functions for the (Numerators of the) Factorial Moments}

Let's define
$$
P(X)=1+ \sum_{i \in S} X^i \quad,
$$
then our functional equation for the grand-generating function, $F(x,y)$ can be written
$$
F(x,y)=xP(F(xy,y) )\quad .
$$
If we want to get generating functions for the first $k$ factorial moments of our random variable $H_n$, we need
the first $k$ coefficients of the Taylor expansion, about $y=1$, of $F(x,y)$. Writing $y=1+z$, and
$$
G(x,z)=F(x,1+z) \quad,
$$
we get the functional equation for $G(x,z)$
$$
G(x,z)=x\, P( G(x+xz,z))  \quad .
\eqno(FE)
$$
Let's write the Taylor expansion of $G(x,z)$ around $z=0$ to order $k$
$$
G(x,z)=\sum_{r=0}^{k} g_r(x) \frac{z^r}{r!} + O(z^{k+1}) \quad.
$$
It follows that
$$
G(x+xz,z)=\sum_{r=0}^{k} g_r(x+xz) \frac{z^r}{r!} + O(z^{k+1}) \quad.
$$

We now do the Taylor expansion of $g_r(x+xz)$ around $x$, getting
$$
g_r(x+xz)=g_r(x) \, + \, g'_r(x)(xz) \, +  \, g''_r(x) \frac{(xz)^2}{2!}  \, + \, \dots \, + \, g_r^{(k)}(x) \frac{(xz)^k}{k!} \, + \, O(z^{k+1}) \quad .
$$ 

Plugging all this into $(FE)$, and comparing coefficients of respective terms of $z^r$ for $r$ from $0$ to $k$
we get $k+1$ extremely complicated equations relating $g^{(j)}_r(x)$ to each other. It is easy to see that
one can express $g_r(x)$ in terms of $g_{s}^{(j)}(x)$ with $s<r$  (and $0 \leq j \leq k$) .

Using implicit differentiation, the derivatives of $g_0(x)$, $g_0^{(j)}(x)$ (where $g_0(x)$ is the same as $f(x)$), can be expressed
as rational functions of $x$ and $g_0(x)$. 
As soon as we get an expression for $g_r(x)$ in terms of $x$ and $g_0(x)$, we can
use calculus to get expressions for the derivatives $g_r^{(j)}(x)$ in terms of $x$ and $g_0(x)$. At the end of the day,
we get expressions for each $g_r(x)$ in terms of $x$ and $g_0(x)$ (alias $f(x)$), and since it is easy to find the first
ten thousand (or whatever) Taylor coefficients of $g_0(x)$, we can get the first ten thousand coefficients of
$g_r(x)$, for all $0 \leq r \leq k$, and get the numerical sequences that will enable us to get the above-mentioned statistical information.

The beauty is that this is all done by the computer! Maple knows calculus.

We can do even better. Using the methods described in [FS], one should be able to get, {\it automatically},
asymptotic formulas for the expectation, variance, and as many moments as desired. Alas, implementing
it in general would have to wait for the future.

For the special case of complete binary trees, everything can be expressed in terms of Catalan numbers, and
hence the asymptotic is easy, and our beloved computer, running the Maple package {\tt TREES.txt} (mentioned above),
obtained the results in the next section.

{\bf Computer-Generated Theorems About the Expectation, Variance, and First Nine Moments for the Total Height on Complete Binary Trees on $n$ Leaves}

See the output file

{\tt http://www.math.rutgers.edu/\~{}zeilberg/tokhniot/oTREES3.txt} .

{\bf Universality}

The computer output, given in the above webpage, proved that for this case, of complete binary trees, the  limits of the first nine scaled moments coincide
{\it exactly} with those found in [EZ2], and given above. 
Confirming, by {\it purely elementary, finitistic methods}, the universality property mentioned above.
We do it for one family at a time, and only for finitely many moments, but on the other hand, we derived
{\it explicit} expressions for the first twelve moments in the case of complete binary trees, and explicit
expressions for the generating functions for other families.

{\bf Conclusion} 

Even more  interesting than the {\it actual} research reported here, it the way that is was obtained. Fully automatically!

\bigskip

{\bf References}

[A1] David Aldous, {\it The continuum random tree II}, {\it The continuum random tree II: an overview},  Stochastic analysis {\bf 167} (1991), 23-70.

[A2] David Aldous, {\it The continuum random tree III}, Ann. Probab. {\bf 21} (1993), 248-289,.

[EZ1] Shalosh B. Ekhad and Doron Zeilberger, 
{\it  Explicit Expressions for the Variance and Higher Moments of the Size of a Simultaneous Core Partition and its Limiting Distribution  },
The Personal Journal of Shalosh B. Ekhad and Doron Zeilberger, \hfill\break
{\tt http://www.math.rutgers.edu/\~{}zeilberg/mamarim/mamarimhtml/stcore.html} \quad .

[EZ2] Shalosh B. Ekhad and Doron Zeilberger, 
{\it   Going Back to Neil Sloane's FIRST LOVE (OEIS Sequence A435): On the Total Heights in Rooted Labeled Trees},
The Personal Journal of Shalosh B. Ekhad and Doron Zeilberger, \hfill\break
{\tt http://www.math.rutgers.edu/\~{}zeilberg/mamarim/mamarimhtml/a435.html} \quad .

[EZ3] Shalosh B. Ekhad and Doron Zeilberger, 
{\it Automatic Proofs of Asymptotic ABNORMALITY (and much more!) of Natural Statistics Defined on Catalan-Counted Combinatorial Families},
The Personal Journal of Shalosh B. Ekhad and Doron Zeilberger, \hfill\break
{\tt http://www.math.rutgers.edu/\~{}zeilberg/mamarim/mamarimhtml/abnormal.html} \quad .

[FS]  Philippe Flajolet Robert Sedgewick, ``{\it  Analytic Combinatorics}'', Cambridge University Press, 2009.
(Free download from Flajolet's homepage.)

[J1] Svante Janson {\it Brownian excursion area, Wright's constants in graph enumeration, and other Brownian areas},
Probability Surveys 3 (2007), 80-145. Available on line from \hfill\break
{\tt https://arxiv.org/abs/	0704.2289}\quad

[J2] Svante Janson {\it Patterns in random permutations avoiding the pattern 132},
Available on line from \hfill\break
{\tt https://arxiv.org/abs/1401.5679}\quad

[J3] Svante Janson {\it Simply generated trees, conditioned Galton–Watson trees, random allocations and condensation},
Probability Surveys 9 (2012), 103-252. Available on line from \hfill\break
{\tt http://www2.math.uu.se/~svante/papers/sj264.pdf} \quad

[J4] Svante Janson {\it The Wiener Index of Simply Generated Trees},
Random Structures and algorithms 22, issue 4 (2003), 337-358. Available on line from \hfill\break
{\tt http://www2.math.uu.se/~svante/papers/sj146.pdf} \quad

[MM] Jean-François Marckert and Abdelkader Mokkadem {\it The depth first processes of Galton--Watson trees converge to the same Brownian excursion},
 The Annals of Probability  Volume 31, number 3 (2003) 1655-1678. Available on line from \hfill\break
{\tt http://projecteuclid.org/euclid.aop/1055425793}\quad .

[RiS] John Riordan and Neil J. A. Sloane, {\it The enumeration of rooted trees by total height},
J. Australian Math. Soc. {\bf 10}(1969),  278-282. Available on line from: \hfill\break 
{\tt http://neilsloane.com/doc/riordan-enum-trees-by-height.pdf} \quad .

[Ro] Dan Romik, ``{\it The Surprising Mathematics of Longest Increasing Subsequences}'', Cambridge University Press,
2015. 

[Z1] Doron Zeilberger, {\it Doron Gepner's Statistics on Words in $\{1,2,3\}$* is (Most Probably) Asymptotically Logistic},
The Personal Journal of Shalosh B. Ekhad and Doron Zeilberger, \hfill\break
{\tt http://www.math.rutgers.edu/\~{}zeilberg/mamarim/mamarimhtml/gepner.html} \quad .

\bigskip
\bigskip
\hrule
\bigskip
Andrew Lohr, Department of Mathematics, Rutgers University (New Brunswick), Hill Center-Busch Campus, 110 Frelinghuysen
Rd., Piscataway, NJ 08854-8019, USA. ajl213 at math dot rutgers dot edu; \quad {\tt http://www.math.rutgers.edu/\~{}ajl213/} \quad .
\bigskip
\hrule
\bigskip
Doron Zeilberger, Department of Mathematics, Rutgers University (New Brunswick), Hill Center-Busch Campus, 110 Frelinghuysen
Rd., Piscataway, NJ 08854-8019, USA. \hfill \break
DoronZeil at gmail dot com  \quad ;  \quad {\tt http://www.math.rutgers.edu/\~{}zeilberg/} \quad .
\bigskip
\hrule

\bigskip
\bigskip

{\bf  First Written:  Feb. 9, 2017; This version: March 6, 2017.}.

\end